\documentclass[12pt]{article} 
\usepackage{amstext,amsfonts,graphics,latexsym}
\newtheorem{thm}{Theorem} \newtheorem{lemma}[thm]{Lemma}
\newtheorem{remark}{Remark} 
 
 \newenvironment{pfof}[1]
{\par\vskip2\parsep\noindent{\sc Proof of\ #1. }}{{\hfill $\Box$}
\par\vskip2\parsep}

\newcommand{\hk}{g_k}
\newcommand{\ko}{K}
\newcommand{\gk}{g_k}
\newcommand{\D}{4}
\newcommand{\E}{2}
\newcommand{\harm}{\frac{4}{k}}
\newcommand{\return}{S_n(0)=\0}
\newcommand{\ratereturn}{|S_n(t)|<n^{.5-1/(\log(n))^{.25+\epsilon}}}
\newcommand{\exc}{\mbox{Exc}}
\newcommand{\ann}{A}
\newcommand{\ret}{R}
\newcommand{\esc}{E}

\newcommand{\zt}{{\mathbb Z}^2}

\newcommand{\prob}{\mbox{\bf P}}
\newcommand{\expe}{\mbox{\bf E}}
\newcommand{\leb}{{\mathcal L}}

\newcommand{\olle}{H\"{a}ggstr\"{o}m}

\newcommand{\Z}{\mathbb Z}
\newcommand{\R}{\mathbb R}
\newcommand{\N}{\mathbb N}
\newcommand{\sk}{k^{10}2^{2k^2}}

\newcommand{\be}{\begin{equation}}
\newcommand{\ee}{\end{equation}}
\newcommand{\0}{{\bf 0}}

\def\qed{\relax\ifmmode\hskip2em \Box\else\unskip\nobreak\hfill
$\Box$\fi}

\newcounter{mycount}
\newenvironment{proof}{\noindent{\sc Proof. }}{{\hfill $\Box$}\par\vskip2\parsep}

\title{Recurrence of Simple Random Walk on $\Z^2$ is Dynamically Sensitive}
\author{Christopher Hoffman}

\begin{document}
\maketitle

\begin{abstract}
Benjamini, H\"{a}ggstr\"{o}m, Peres and Steif \cite{ds}
introduced the concept of a dynamical random walk.  This is a
continuous family of random walks, $\{S_n(t)\}_{n \in \N, t \in
\R}$. Benjamini et.\ al. proved that if $d=3$ or $d=4$ then there
is an exceptional set of $t$ such that $\{S_n(t)\}_{n \in \N}$
returns to the origin infinitely often. In this paper we consider
a dynamical random walk on $\zt$. We show that with probability
one there exists $t\in \R$ such that $\{S_n(t)\}_{n \in \N}$ never
returns to the origin.
This exceptional set of times has dimension one.
This proves a conjecture of Benjamini et.\ al. \cite{ds}.
\end{abstract}

\section{Introduction}
We consider a dynamical simple random walk on $\zt$.  Associated
with each $n$ is Poisson clock.  When the clock rings the $n$th
step of the random walk is replaced by an independent random
variable with the same distribution. Thus for any fixed $t$ the
distribution of the walks at time $t$ is that of simple random
walk on $\zt$ and is almost surely recurrent.

We prove that with probability one there exists a (random) set of
times $t$ such that $S_n(t) \neq \0 \ \forall \ n\in \N$. Thus we say that
recurrence of simple random walk on $\zt$ is {\bf dynamically
sensitive}.

More formally let $\{Y_n^m\}_{m,n \in \N}$ 
be uniformly distributed i.i.d.\ random
variables chosen from the set $\{(0,1),(0,-1),(1,0),(-1,0)\}$. Let
$\{\tau_n^{(m)}\}_{m \geq 0,n \in \N}$ be an independent Poisson process of
rate one and $\tau_n^{(0)}=0$ for each $n$. Define
$$X_n(t)=Y_n^m$$
for all $t \in [\tau_n^{(m)},\tau_n^{(m+1)}).$ Let
$$S_n(t)=\sum_{i=1}^{n}X_i(t).$$
Thus for each $t$ the random variables $\{X_n(t)\}_{n\in \N}$ are
i.i.d.

Define the exceptional set of times
$$\exc =\{t:\     S_n(t) \neq \0 \forall \ n\}.$$

Our main result is
\begin{thm} \label{main}
 $$\prob(\exc \neq \emptyset)=1.$$
Moreover,
$\exc$ has dimension 1 a.s.
\end{thm}

\begin{remark} \label{remark1}
Our methods can be used to 
calculate a rate of escape.  For any 
$\alpha<1/2$ there is a set of $t$ such that $|S_n(t)|>n^{\alpha}$
for all $n$.
The limits of are method yield that with probability one there
is a time $t$ such where the rate of escape is at least
$$|S_n(t)|>n^{.5-1/(\log n)^{1/4 + \epsilon}}$$ for all $n$.

\end{remark}

Benjamini, \olle, Peres and Steif introduced the concept of dynamical 
random walk and showed that the strong law of
large numbers and the law of iterated logarithms are satisfied for
all times almost surely \cite{ds}.  Thus these properties are said to be {\bf
dynamically stable}. They also proved that in dimensions 3 and 4
that the transience of simple random walk is dynamically sensitive
and in dimensions 5 and higher that transience is dynamically
stable.  Levin, Khoshnevesian and Mendez have studied other properties 
of dynamical random walks \cite{levinone} and 
\cite{levintwo}.
\olle, Peres and Steif studied similar questions of dynamic 
stability and sensitivity for percolation \cite{dp}.

Dynamical random walk and the results in this paper are 
related to several other topics in probability.
Most closely related to the work in this paper is a result of 
Adelman, Burdzy and Pemantle about
sets missed by three dimensional Brownian motion \cite{abp}.  
The projection of Brownian motion on $\R^3$ onto a fixed plane yields
Brownian motion in the plane which is neighborhood recurrent. For
a fixed plane the projection of almost every Brownian path onto the plane is
neighborhood recurrent.  They proved that with probability one
there is a (random) set of exceptional planes such that the set of
times that the projected path is in any bounded set is bounded.

The questions studied about dynamical random walks and dynamical percolation
have a strong resemblence to questions of quasi-everywhere properties of 
Brownian paths.  These are properties that hold simultaneously 
for every cross section of a Brownian sheet with probability one.  See
\cite{BS1} and \cite{BS2}.

\section{Outline}
We start by introducing some notation.
Let $s_0=1$ and $s_k=\sk$ for $k\geq 1$.  
This is a sequence of stopping times. Define
the event $\ret_k(t)$ to be
$$\ret_k(t)=\{\exists n \in \{s_{k-1},\dots , s_k\} \text{ such that } 
\return\}.$$
%
For $x \in \zt$ we use the standard notation
$|x|=\sqrt{(x_1)^2+(x_2)^2}$. Define the annulus
$$\ann_k=\{x \in\zt :\ 2^{k^2} \leq |x| \leq k^{10}2^{k^2}\}.$$
Define the event $G_k(t)$ to be
$$G_k(t)=\{S_{s_k}(t) \in \ann_k \}.$$
Also define the events
 $G_k(0,t)= G_k(0) \cap G_k(t)$ and
 $\ret_k(0,t)= \ret_k(0) \cap \ret_k(t)$.

$$\esc_M(0)= \left(\cap_{1}^{M}G_k(0)\right) 
		\setminus \left(\cup_1^{M} \ret_k(0)\right).$$
$$\esc_M(0,t)= \esc_M(0) \cap \esc_M(t).$$
We will show in Lemma \ref{smallk}
that there is an integrable function $f(t)$ such that
for all $M$
\begin{equation} \label{thisbound}
\frac{\prob(\esc_M(0,t))}{(\prob(\esc_M(0)))^2}<f(t).
\end{equation}
Theorem \ref{main} follows from Lemma \ref{smallk} by the
second moment method.

We obtain (\ref{thisbound}) by mulitplying together conditional
probabilities. Some of our bounds will hold only when $k$ is
sufficiently large compared with $1/t$.  For this reason we define
$\ko=\ko(t)$ to be the unique integer such that
\begin{equation} \label{tkrelation}
1+|\log t| > \ko \geq |\log t|
\end{equation}
for $t<1$ and $K=0$ if $t \geq 1$.
Our main lemma is the following.

\begin{lemma} \label{summary}
There is a positive sequence $\gk$ such that 
\begin{enumerate}
\item $\sum_1^\infty \gk< \infty$, 
\item $(\prob(\esc_k(0)|\esc_{k-1}(0)))^2>1-\harm-\hk$
	for all $k>1$, and \label{summary2}
\item $\prob(\esc_k(0,t)|\esc_{k-1}(0,t))<1-\harm +\gk$ for all $k>\ko$.
	\label{summary1}
\end{enumerate}
\end{lemma}
The next section is dedicated to proving Lemma \ref{summary}.  In
the last section we show how Lemma \ref{summary} implies Theorem
\ref{main}. 

We end this section with a few notes about notation.
We use $\log (n) =\log_{2}(n)$.  We use $C$ as a
generic constant whose value may increase from line to line and
lemma to lemma. In many of the proofs in the next section we use
bounds that only hold for sufficiently large $k$.  This causes no
problem since it will be clear that 
we can always choose $C$ such that the lemma is true
for all $k$.

\section{Proof of Lemma \ref{summary}}
For the rest of the paper we will use
the following notation for conditional probabilities.
Let
$$\prob^{x,k-1}(*)=\prob\left(* |\
      S_{s_{k-1}}(0)=x \right) $$ and
$$\prob^{x,y,k-1}(*)=\prob\left(* |\
      S_{s_{k-1}}(0)=x \mbox{ and } S_{s_{k-1}}(t)=y\right) $$
The two main parts of the proof of Lemma \ref{summary} are Lemma \ref{mainlem} where we get upper and lower bounds
on $\prob^{x,k-1}(\ret_k(0))$ and Lemma
\ref{bothreturn} where we get an upper bound on
 $\prob^{x,y,k-1}(\ret_k(0,t)) $.
The main tool that we use are bounds on the probability that
simple random walk started at $x$ returns to the
origin before exiting the ball of radius $n$ and center at the
origin. The probability of this is calculated in Proposition 1.6.7 on
page 40 of \cite{lawler}. We use only a weak version of the result
there.

Let $\eta$ be the  smallest $m>0$ such that $S_m(0)=\0$ or
$|S_m(0)| \geq n$.
\begin{lemma} \label{lawler}
 There exists $C$ such that for all $x$ with $0<|x|<n$
$$\frac{\log(n)-\log|x|-C}{\log(n)}\leq \prob(S_\eta(0)=\0|S_0(0)=x)\leq
\frac{\log(n)-\log|x|+C}{\log(n)}.$$
\end{lemma}

We will frequently use the following standard bounds.  
\begin{lemma} \label{leave}
There exists $C$ such that for all $x \in \zt$, $n \in \N$ and $m<\sqrt{n}$
\begin{equation} \label{le}
\prob(\exists n'<n:S_{n'}(0)>m\sqrt{n})\leq \frac{C}{m^2}
\end{equation} and
$$\prob(|S_{n}(0)-x|<\frac{\sqrt{n}}{m})\leq \frac{C}{m^2}.$$
\end{lemma}

\begin{proof}
If  $|S_{n'}(0)|>m\sqrt{n}$
then one component of the random walk has absolute value bigger than
$m\sqrt{n}/2$ .
Thus the left hand side of (\ref{le})
is at most than four times the probability that one dimensional 
simple random walk is ever more than $m\sqrt{n}/2$ 
away from the origin during the first $n$
steps. The probability that a one dimensional 
simple random walk has ever been larger
than $m\sqrt{n}/2$ in the first $n$ steps is at most twice 
the probability that one dimensional 
simple random walk is greater than $m\sqrt{n}/2$ 
after $n$
steps.  Chebyshev's inequality then gives the first bound.

To bound the probability that $|S_{n}(0)-x|$ is too small we note
that since $m<\sqrt{n}$ 
the number of $y \in \zt$ such that  $|y-x|<\frac{\sqrt{n}}{m}$  
is less
than $\frac{10n}{m^2}$.  There is $C$ such that 
for any $n \in \N$ and $z \in \zt$ the
probability that $S_{n}(0)=z$ is less than
$C/n$. 
\end{proof}

\begin{lemma} \label{mainlem}
There exists $C$ such that for any $k$ and any $x \in \ann_{k-1}$
$$\frac{\E}{k} -\frac {C\log k}{k^2}\leq \prob^{x,k-1} (\ret_k(0)) 
	\leq \frac{\E}{k} +\frac {C\log k}{k^2}.$$
\end{lemma}

\begin{proof}
If the random walk returns to $\0$ in less than $s_k$ steps then
either it returns to $\0$ before exiting the ball of radius
$\sqrt{s_k} \log(s_k)$ or it exits the ball in less than $s_k$
steps.  Thus by Lemmas \ref{lawler} and \ref{leave} our
upper bound is
\begin{eqnarray*}
 &<& \frac{\log(\sqrt{s_k}\log(s_k))-\log(2^{(k-1)^2})+C}
		{\log(\sqrt{s_k}\log(s_k))}
+\frac{C}{(\log s_k)^2}\\
 &<&\frac{5\log k+k^2+C \log k -(k-1)^2+C}
	{5\log k+k^2 +\log (2k^2+10\log k)}
 +\frac{C}{k^4}\\
 &<&\frac{\E k+C\log k}{k^2}.
 \end{eqnarray*}
If the random walk returns to $\0$  after $s_{k-1}$ but 
before exiting the ball of
radius $\sqrt{s_k}/ \log(s_k)$ and it is outside the ball of
radius $\sqrt{s_k}/ \log(s_k)$ at time $s_k$ then it has
returned to $\0$ between times $s_{k-1}$ and $s_k$. Thus by Lemmas
\ref{lawler} and \ref{leave} our lower bound is
\begin{eqnarray*}
&>&\frac{\log(\sqrt{s_k}/\log(s_k))-\log((k-1)^{10}2^{(k-1)^2})-C}
	{\log(\sqrt{s_k}/\log(s_k))}
-\frac{C}{(\log(s_k))^2}\\
 &>&\frac{5\log k+k^2-C\log k-10\log(k-1)-(k-1)^2-C}
	{5\log k+k^2 -\log \log s_k}
-\frac{C}{k^4}\\
 &>&\frac{\E k-C\log k}{k^2}.
 \end{eqnarray*}

\end{proof}

\begin{lemma} \label{gbound} For any $k$ and $x \in \ann_{k-1}$
$$\prob^{x,k-1} ((G_k(0))^C) \leq \frac{C}{k^{10}}.$$
\end{lemma}

\begin{proof}
This follows directly from Lemma \ref{leave}
\end{proof}

Now we start to bound the probability that both walks return to the origin
between times $s_{k-1}$ and $s_k$.
We first need the following lemma.

\begin{lemma} \label{returnafter}
There exists $C$ such that for all $k$, $n\geq
s_{k-1}+s_k/2^{10k}$, for all $I \subset \{1,\dots,n\}$ with
$|I|\geq s_k/2^{10k}$ and for all
 $\{x_i(t)\}_{i \in \{1,\dots,n\} \setminus I}$
$$    \prob \left(\exists j \in \{n,\dots,s_k\}
    \mbox{ such that } S_j(t)=\0\right |\
    \{x_i(t)\}_{i \in \{1,\dots,n\}\setminus I} )\leq
    \frac{C}{k}.$$
\end{lemma}

\begin{proof}
Rearranging the first $n$ steps of a random walk does not change
the random walk after time $n$. The probability is largest when
$n$ and $|I|$ are as small as possible. Thus it causes no loss of
generality to assume that $n=s_{k-1}+s_k/2^{10k}$ and
$I=\{s_{k-1}+1,\dots,n\}$.

If the event happens then either
\begin{enumerate}
\item $|S_{n}(t)| \leq \frac{\sqrt{|I|}}{\log{|I|}}$,
\item $|S_{n}(t)| > \frac{\sqrt{|I|}}{\log{|I|}}$ and
$$\inf \{j: j>n \mbox{ and }  S_j(t)=\0\} <
  \inf \{j: j>n \mbox{ and } |S_j(t)|>\sqrt{s_k}\log(s_k)\}$$
 or
\item there exists $j'$ such that $n<j' <s_k$
    and $|S_{j'}(t)-S_n(t)| \geq .5\sqrt{s_k}\log{s_k}$.
\end{enumerate}

The probability of the first and third events are bounded
by Lemma \ref{leave}. The probability of the second event 
is bounded by Lemma \ref{lawler}. Thus our upper bound is
\begin{eqnarray*}
 &<& \frac{C}{(\log |I|)^2}+
 \frac{\log(\sqrt{s_k}\log(s_k))-\log(\frac{\sqrt{|I|}}{\log |I|})+C}{\log(\sqrt{s_k}\log(s_k))}
+\frac{C}{(\log s_k)^2}\\
 &<&\frac{C}{(\log |I|)^2}+
  \frac{.5\log(s_k)+\log(\log s_k)-(.5\log s_k - \log(2^{5k})-\log(\log |I|)) +C}
	{.5\log(s_k)+\log(\log(s_k))}\\
 &<&\frac{C}{(2k^2-10k)^2}+
  \frac{Ck +C\log k+C}{k^2}\\
 &<&\frac{C}{(2k^2-10k)^2}+
  \frac{C k}{k^2}\\
 &<&\frac{C}{k}.
 \end{eqnarray*}

\end{proof}

\begin{lemma}  \label{bothreturn}
There exists $C$ such that for any $t$, any $k>\ko(t)$ and any $x,y
\in \ann_{k-1}$
$$\prob^{x,y,k-1} (\ret_k(0,t))
    \leq \frac{C}{k^2}.$$
\end{lemma}

\begin{proof}
Let
$$I\subset \left\{s_{k-1},\dots,s_{k-1}+2^{2(k-1)^2}/k^2\right\}$$
be the set of $i$ such that conditioned on the Poisson process,
$X_i(t)$ and $X_i(0)$ are independent. Let $B$ be the event that
there exists $n$ such that
$$s_{k-1} \leq n\leq s_{k-1}+2^{2(k-1)^2}/k^2$$
such that
$S_n(0)=\0$. Let $D$ be the event that there exist $n$ and $n'$
such that
\begin{enumerate}
\item $s_{k-1}+2^{2(k-1)^2}/k^2< n \leq n'\leq s_k$
\item $S_n(0)=S_{n'}(t)=\0$ and
\item $|I| \geq s_k/2^{10k}$.
\end{enumerate}

If $\ret_k(0,t)$ occurs then either the first return happens before step 
$s_{k-1}+2^{2(k-1)^2}/k^2$ or after that step.  
The probability that the first return is before
is bounded by twice the probability of $B$.  If the first reutrn is after
then either $|I|< s_k/2^{10k}$ or $|I|\geq  s_k/2^{10k}$. The 
probability that the first return is after 
and $|I|$ is large is bounded by twice the probability of $D$. 
Thus we get that
\begin{eqnarray*}
\prob^{x,y,k-1} (\ret_k(0,t))
    &\leq &2\prob^{x,y,k-1}(B)+\prob^{x,y,k-1}(|I|<s_k/2^{10k})
        +2\prob^{x,y,k-1}(D).
\end{eqnarray*}
By (\ref{tkrelation}) $\min(1,t)\geq 1/2^{K}.$ As $k>K$
this implies the expected size of $|I|$ is
\begin{equation} \label{star4}
(1-e^{-t})\frac{2^{2(k-1)^2}}{k^2}
    > \frac12 \min(1,t)\frac{2^{2(k-1)^2}}{k^2}
    \geq \frac{2^{2(k-1)^2}}{2^{\ko +1}k^2}
    > \frac{2^{2k^2}}{2^{6k}}
    > 2\frac{s_k}{2^{10k}}.
\end{equation}
 Thus the probability
that $|I|<s_k/2^{10k}$ is at most $C/k^2$.

By Lemma \ref{leave} the conditional probability of $B$ is 
bounded by 
$C/k^2$.
In order for $D$ to happen we first need that the event
$\ret_k(0)$ occurs.  By Lemma \ref{mainlem} the probability of
this is bounded by $C/k.$ Now we condition on the following events
\begin{enumerate}
\item $\{X_i(0)\}_{i\geq 0}$
\item the Poisson process,
\item $|I|\geq 2^{-10k}s_k$, and
\item $\{X_i(t)\}_{i \in \{1,\dots,n\}\setminus I}$
\end{enumerate}
 and bound the probability that there exists $n' \in
\{n, \dots ,s_k\}$ with $S_{n'}(t)=\0$.   

By the first condition in the definition of $D$ and  (\ref{star4}) 
$$ n>s_{k-1}+\frac{2^{(k-1)^2}}{k^2}>s_{k-1}+s_k/2^{10k}$$
and 
Lemma \ref{returnafter} applies.  Thus the
conditional probability of $D$ given $\ret_k(0)$ is at most $C/k$.

Putting this together we get
\begin{eqnarray*}
\prob^{x,y,k-1} (\ret_k(0,t))
    &\leq &   2\prob^{x,y,k-1}(B)+\prob^{x,y,k-1}(|I|<\frac{s_k}{2^{10k}})
	+ 2\prob^{x,y,k-1}(D)   \\
    &\leq & \frac{C}{k^2}+\frac{C}{k^2}+2\left(\frac{C}{k} \right)
    \left(\frac{C}{k} \right) \\
    &\leq & \frac{C}{k^2}.
\end{eqnarray*}
\end{proof}

\begin{pfof}{Lemma \ref{summary}}
We let 
$$\gk=\frac{C\log k}{k^2}.$$  
Clearly this satisfies the summability condition.
Note that if $\esc_{k-1}(0)$ occurs then $G_{k-1}(0)$ occurs and
$S_{s_{k-1}}(0)\in \ann_{k-1}$.  Since simple random walk is Markovian,
for any $x \in \ann_{k-1}$
\begin{eqnarray*}
\prob(\esc_k(0)|\esc_{k-1}(0)\mbox{ and }S_{s_{k-1}}(0)=x)
 &=& \prob^{x,k-1}((\ret_k(0))^C \cap G_k(0)).
\end{eqnarray*}
Along with Lemmas \ref{mainlem} 
and  \ref{gbound} 
this tells us that
\begin{eqnarray*}
\prob(\esc_k(0)|\esc_{k-1}(0))
 &\geq&  \min_{x \in \ann_{k-1}}\prob
	 \left(\esc_{k}(0)|\esc_{k-1}(0)\mbox{ and }S_{s_{k-1}}(0)=x\right)\\
 &\geq&  \min_{x \in \ann_{k-1}}\prob^{x,k-1}((\ret_k(0))^C \cap G_k(0))\\
 &\geq&  \min_{x \in \ann_{k-1}}\prob^{x,k-1}((\ret_k(0))^C )
         	-\max_{x \in \ann_{k-1}}\prob^{x,k-1} \left((G_k(0))^C\right)\\
 &\geq&  1-\max_{x \in \ann_{k-1}}\prob^{x,k-1}\left(\ret_k(0)\right)
		-\frac{C}{k^{10}}\\
 & \geq& 1-\left( \frac{2}{k} +\frac{C\log k}{k^2}\right) -\frac{C}{k^{10}}\\
 & \geq& 1-\frac{2}{k} -\frac{C\log k}{k^2}.
\end{eqnarray*}
Squaring both sides yields condition \ref{summary2} of Lemma \ref{summary}.

Also note that if $\esc_{k-1}(0,t)$ occurs then $G_{k-1}(0,t)$ occurs and
$$S_{s_{k-1}}(0),S_{s_{k-1}}(t) \in \ann_{k-1}.$$  
Since dynamic random walk is Markovian, for any $x,y \in \ann_{k-1}$
 
\begin{eqnarray*}
\lefteqn{
\prob(\esc_k(0,t)|\esc_{k-1}(0,t)\mbox{ and }
	S_{s_{k-1}}(0)=x, S_{s_{k-1}}(t)=y)}&&\\
 &=& \prob^{x,y,k-1}((\ret_k(0))^C \cap (\ret_k(t))^C \cap G_k(0,t)).
\end{eqnarray*}
Combining this with 
 Lemmas \ref{mainlem}
and \ref{bothreturn}
we get that
\pagebreak[0]	
\begin{eqnarray*}
\lefteqn{\prob(\esc_k(0,t)|\esc_{k-1}(0,t))}&&\\ 
 &\leq &  \max_{x,y \in \ann_{k-1}}\prob^{x,y,k-1}
		( (\ret_k(0))^C \cap (\ret_k(t))^C \cap G_k(0,t))\\
  &\leq& \max_{x,y \in \ann_{k-1}}\prob^{x,y,k-1}\left((\ret_k(0))^C
    \cap (\ret_k(t))^C \right)\\
 &\leq& 1-2\min_{x \in \ann_{k-1}}\prob^{x,k-1}(\ret_k(0))
+\max_{x,y\in \ann_{k-1}}\prob^{x,y,k-1}(\ret_k(0,t))\\
& \leq& 1-2 \left(\frac{2}{k} -\frac{C\log k}{k^2} \right) +\frac{C}{k^2}\\
& \leq& 1- \frac{4}{k} -\frac{C\log k}{k^2}.
\end{eqnarray*}
This proves condition \ref{summary1} of Lemma \ref{summary}.  
\end{pfof}

\section{Proof of Theorem \ref{main}}

Define
$$f(t,M) =\frac{\prob(\esc_M(0,t))}
            {(\prob (\esc_M(0)))^2}$$

\begin{lemma} \label{smallk}
There exists $C$ such that for any $t$ and any $M$
\begin{equation} \label{star3}
f(t,M) <C(1+|\log t|)^{\D}.
\end{equation}
\end{lemma}

\begin{proof}
 Choose $n$ such that
 $$\harm +g(k)<.5$$
for all $k\geq n$. By Lemma \ref{summary} 
\begin{equation} \label{k_0}
 f(t,M) \leq \frac{1}{(\prob(\esc_n(0)))^2}
 \prod_{k=n+1}^{\ko}\frac{1}{1-\harm-\hk}
 \prod_{\ko+1}^{M}\frac{1-\harm+\gk}{1-\harm-\hk}.
 \end{equation}
 The inequality
$$-x^2-x<\ln(1-x)<-x$$
holds for all $x \in (0,.5)$. Thus
\begin{eqnarray}
\ln \left( \prod_{k=n+1}^{\ko}\frac{1}{1-\harm -\hk}\right)
    & = &-\sum_{k=n+1}^{\ko}\ln \left(1-\harm -\hk\right)\nonumber\\
    & < &\sum_{k=n+1}^{\ko}\harm +\hk+ \left(\harm +\hk \right)^2\nonumber\\
    & < &C+\sum_{k=n+1}^{\ko}\harm \nonumber\\
    & < &C+\D \ln(\ko). \nonumber
\end{eqnarray}
By exponentiating both sides and (\ref{tkrelation}) we get
\begin{equation} \label{line1}
 \prod_{k=n+1}^{\ko}\frac{1}{1-\harm -\hk}
 \leq C\ko^{\D}
 \leq C(1+|\log t|)^{\D}.
\end{equation}

\begin{eqnarray*}
\ln\left(\prod_{\ko+1}^{M} \frac{1-\harm +\gk}{1-\harm-\hk}\right)
 &  \leq & \sum_{\ko+1}^{\infty}\ln(1-\harm+\gk)
    -\sum_{\ko+1}^{\infty} \ln (1-\harm -\hk)\\
 &  \leq& \sum_{\ko+1}^{\infty} -\harm+\gk
    - \left(-\left(\harm +\hk\right) -\left(\harm +\hk\right)^2 \right)\\
 &  \leq & \sum_{\ko+1}^{\infty} 2\gk +\frac{16}{k^2}
    +\frac{8\hk}{k} +\hk^2\\
 & \leq & C.
\end{eqnarray*}
Exponentiating both sides we get for all $M$
\begin{equation} \label{line2}
 \prod_{\ko+1}^{M} \frac{1-\harm +\gk}{1-\harm-\hk}\leq  C.
\end{equation}
 Putting together  (\ref{k_0}), (\ref{line1}) and
(\ref{line2}) we get
$$f(t,M) \leq \frac{1}{(\prob(\esc_n(0)))^2}C(1+|\log t|)^{\D} C
   \leq    C(1+|\log t|)^{\D}.$$
\end{proof}

\begin{pfof}{Theorem \ref{main}}
Define 
$$T_M=\{t:t \in [0,1] \mbox{ and }\esc_M(t) \mbox{ occurs}\}$$ 
and 
$$T=\cap_1^{\infty}\overline{T_M}.$$
Now we show that $T$ is contained in the union of $\exc$ and 
the countable set $$\Lambda=(\cup_{n,m}\tau_n^{(m)})\cup 1.$$  
If $t \in \cap_1^{\infty} T_M$ then $t \in \exc$.
So if $t \in T \setminus \exc$ then $t$ is contained in the boundary of $T_M$ 
for some $M$. For any $M$ the boundary of $T_M$ is contained in $\Lambda$. 
Thus if $t \in T \setminus \exc$ then $t\in \Lambda$ and
$$T \subset \exc \cup \Lambda.$$ 
As $\Lambda$ is countable 
if $T$ has dimension one with positive probability then so does $\exc$.

By Lemma \ref{smallk} there exists $f(t)$ such that
 $$\int_0^1 f(t) dt < \infty$$ and for all $M$
 \begin{equation} \label{downlow}
 \frac{\prob(\esc_M(0,t))}{(\prob(\esc_M(0)))^2}<f(M,t)< f(t).
 \end{equation}
Let $\leb(*)$ denote Lebesgue measure on $[0,1]$.  Then we get
\begin{eqnarray}
\expe(\leb(T_M)^2)
 & =    & \int_0^1 \int_0^1 \prob(\esc_M(r,s)) dr \times ds \label{fubini}\\
 & \leq & \int_0^1 \int_0^1 \prob(\esc_M(0,|s-r|)) dr \times ds \label{inv}\\
 & \leq & \int_0^1 2\int_0^1 \prob(\esc_M(0,t)) dr \times dt \nonumber\\
 &\leq  & 2 \int_0^1 f(t)  \prob(\esc_M(0))^2 dt \label{downlower}\\
 &\leq  & 2  \prob(\esc_M(0))^2\int_0^1 f(t)  dt. \nonumber
\end{eqnarray}
The equality  (\ref{fubini}) is true by Fubini's theorem,
(\ref{inv}) is true because 
$$\esc_M(a,b)=\esc_M(b,a)=\esc_M(0,|b-a|)$$
and (\ref{downlower}) follows from (\ref{downlow}).

Then
\begin{eqnarray}
2  \prob(\esc_M(0))^2\int_0^1 f(t)  dt & \geq &
\expe(\leb(T_M)^2) \nonumber\\
 &\geq& \left(\frac{\expe(\leb(T_M))}{\prob(T_M \neq \emptyset)}\right)^2
    \cdot \prob(T_M \neq \emptyset)\nonumber\\  
 & \geq & \frac{1}{\prob(T_M \neq \emptyset)}\left(
    \expe(\leb(T_M))\right)^2 \nonumber\\
 & \geq & \frac{1}{\prob(T_M \neq \emptyset)}\prob(\esc_M(0))^2. \nonumber
\end{eqnarray}
Thus for all $M$
$$\prob(T_M \neq \emptyset)
	\geq \frac{1}{2\int_{0}^{1}f(t) \ dt}>0. $$
As $T$ is the intersection of the nested sequence of compact sets $T_M$
$$
\prob(T \neq \emptyset) 
	= \lim_{M \to \infty} \prob(T_M  \neq \emptyset) \geq 
    \frac{1}{2\int_{0}^{1}f(t) \ dt }.$$

Now we show that the dimensions of $T$ and $\exc$ are one.
By Lemma 5.1 of \cite{peres} for any $\beta<1$
there exists a random nested sequence of 
compact sets $F_{k} \subset [0,1]$ 
such that 
\begin{equation} \label{star2}
\prob(r \in F_{k})\geq C(s_k)^{-\beta}
\end{equation}
and
\begin{equation}\label{star1}
\prob(r,t \in F_{k})\leq C(s_k)^{-2\beta}|r-t|^{-\beta}.
\end{equation}
These sets also have the property that for any set $T$ if
\begin{equation} \label{yuval}
\prob(T \cap (\cap_1^{\infty} F_k) \neq \emptyset)>0
\end{equation}
then $T$ has dimension at least $\beta$.
We construct $F_k$ to be independent of the dynamical random walk.
So by (\ref{star3}), (\ref{star2}) and  (\ref{star1}) we get
\begin{equation} \label{dim}
\frac{\prob(r,t \in T_M \cap F_M)}{\prob(r \in T_M \cap F_M)^2}
	\leq C(1+\big|\log |r-t|\big|)^{\D}|r-t|^{-\beta}.
\end{equation}

The same second moment argument 
as above and (\ref{dim}) implies that with positive probability
$T$ satisfies
(\ref{yuval}). 
Thus $T$ has dimension $\beta$ with positive probability.
As 
$$T \subset (\exc \cap [0,1]) \cup \Lambda,$$
and $\Lambda$ is countable,
the dimension of the set of $\exc \cap [0,1]$ is at least $\beta$ 
with positive probability.
By the ergodic theorem 
the dimension of the set of $\exc$ is at least $\beta$ with probability one.
As this holds for all $\beta<1$ the dimension of $\exc$ is one a.s.
\end{pfof}

Finally we briefly state how to modify the proof to calculate  
the rate of escape mentioned in Remark \ref{remark1}.  
For any $\epsilon>0$ 
we replace the event $\ret_k(t)$ with 
$$\ret^{\epsilon}_k(t)=\left\{\exists n \in \{s_{k-1},\dots , s_k\} 
\text{ such that }  \ratereturn \right \}.$$
Instead of Lemma \ref{lawler} we use Exercise 1.6.8 of \cite{lawler}.  
The proof goes through with only minor modifications.

\section*{Acknowledgments}
I would like to thank David Levin and Yuval Peres for introducing
me to this problem and for useful conversations.


\begin{thebibliography}{AAA}

\bibitem{abp}
Adelman, O., Burdzy, K. and Pemantle, R. (1998). Sets avoided by
Brownian motion.  Ann. Probab. 26 429--464.

\bibitem{ds}
Benjamini, Itai; H\"{a}ggstr\"{o}m, Olle; Peres, Yuval; Steif,
Jeffrey E. Which properties of a random sequence are dynamically
sensitive? Ann. Probab.  31  (2003),  no. 1, 1--34.

\bibitem{BS1}
Fukushima, Masatoshi.
Basic properties of Brownian motion and a capacity on the Wiener space.
J. Math. Soc. Japan 36 (1984), no. 1, 161--176.

\bibitem{dp}
H\"{a}ggstr\"{o}m, O., Peres, Y. and Steif, J. E. (1997).
Dynamical percolation.  Ann. Inst. H. Poincar\'e Probab. Statist.
33 497--528.

\bibitem{lawler} Lawler, Gregory F. {\em Intersections of random walks.
Probability and its Applications.} Birkhäuser Boston, Inc.,
Boston, MA, 1991.

\bibitem{levinone}
Levin, D., Khoshnevesian D. and Mendez, P. Exceptional Times and
Invariance for Dynamical Random Walks. math.PR/0409479 
to appear in 
Probability Theory and Related Fields.

\bibitem{levintwo}
Levin, D., Khoshnevesian D. and Mendez, P. 
On Dynamical Gaussian Random Walks.  math.PR/0307346
to appear in Annals of Probability.

\bibitem{BS2}
Penrose, M. D.
On the existence of self-intersections for quasi-every Brownian path in space.
Ann. Probab. 17 (1989), no. 2, 482--502.

\bibitem{peres}
Peres, Yuval. 
Intersection-equivalence of Brownian paths and certain branching processes.  
Comm. Math. Phys.  177  (1996),  no. 2, 417--434.


\end{thebibliography}
\end{document}